\documentclass{amsproc}

\usepackage{a4wide}
\usepackage{amsmath}
\usepackage{enumerate}

\usepackage{amsfonts, amssymb,amsthm,amsgen}

\def\ker{\mathop{{\rm ker}}\nolimits}
\def\spa{\mathop{{\rm span}}\nolimits}
\def\pr{\mathop{{\rm pr}}\nolimits}
\def\Real{\mathbb{R}}

\def\be{\begin{equation}}
\def\ee{\end{equation}}

\newtheorem{lem}{Lemma}
\newtheorem{prop}{Proposition}

\title{On the de~Rham-Wu decomposition for Riemannian and Lorentzian manifolds}

\begin{document}

\author{Anton S. Galaev}

\maketitle

\begin{abstract} It is explained how to find the de~Rham
decomposition of a Riemannian manifold and the Wu decomposition of
a Lorentzian manifold. For that it is enough to find parallel
symmetric bilinear forms on the manifold, and do some  linear
algebra. This result will allow to compute the connected holonomy
group of an arbitrary Riemannian or Lorentzian manifold.



\end{abstract}

\section{Introduction}

The classification of connected holonomy groups  of indecomposable
Riemannian manifolds is a classical result that has many
applications both in geometry and theoretical physics, in
particular, in string theory compactifications and M-theory, see
\cite{Besse,Cecotti,Gubser,Joyce07,McIn2,Pa} and references
therein.

The classification of connected holonomy groups  of indecomposable
Lorentzian manifolds is available as well \cite{ESI,H-L}. Holonomy
groups of 4-dimensional Lorentzian manifolds and their relation to
General Relativity are studied e.g. in \cite{H-L,JJ}. Recently an
attention to the holonomy groups of Lorentzian manifolds of
arbitrary dimension is paid in the physical literature
\cite{BCH,BCH1,Cecotti,G-P,FF00}.

It is important to find the connected holonomy group of an
arbitrary Riemannian or Lorentzian manifold. De~Rham and Wu
theorems allow to decompose (at least locally) a Riemannian or
Lorentzian manifold into a product of indecomposable manifolds.
This implies that the holonomy group of the manifold is a product
of the holonomy groups of indecomposable manifolds. Thus in order
to find the holonomy group of an arbitrary Riemannian or
Lorentzian manifold one should first to find its de~Rham or Wu
decomposition. The proof of the de~Rham and Wu decomposition
theorems assumes that the holonomy group is known, so a priori it
is unclear how to find these decompositions without knowing the
holonomy group.

In this paper we give algorithms that allow to find the de~Rham
decomposition of a Riemannian manifold and the Wu decomposition of
a Lorentzian manifold. In order to find these decompositions, it
is enough to find parallel symmetric bilinear forms on the
manifold, and then do some  linear algebra in the tangent space at
a point of the manifold. Consequently, the algorithms can be
computerized, e.g. as a part of the package DifferentialGeometry
for Maple \cite{Anderson}.

In \cite{Besse}, it is explained how to find the connected
holonomy group of an indecomposable Riemannian manifold, for that
one can analyze parallel differential forms on the manifold. In
another paper we will explain how to find the holonomy group of
an indecomposable Lorentzian manifold. This and the results of the
present paper will provide the complete algorithm that allows to
find the connected holonomy group of an arbitrary Riemannian or
Lorentzian manifold.

{\bf Acknowledgements.} I am grateful  to I.~M.~Anderson for
taking my attention to the problem of finding the algorithm for
computing the holonomy group of a Lorentzian manifold and the
tutorials on Maple.  I am thankful to D.V.~Alekseevsky  for
helpful suggestions.

\section{Holonomy groups}

The theory of holonomy groups of pseudo-Riemannian manifolds can
be found e.g. in \cite{Besse,Joyce07}.

Let $(M,g)$ be a connected pseudo-Riemannian manifold of signature
$(r,s)$ ($r$ is the number of minuses in the signature of the
metric $g$). We will be interested in the case of Riemannian
manifolds ($r=0$, i.e. $g$ is positive definite) and in the case
of Lorentzian manifolds ($r=1$).

Denote by $\nabla$ the Levi-Civita connection on $M$ defined by
the metric $g$; $\nabla$ is the unique torsion-free linear
connection on $M$ that preserves the metric $g$: $\nabla g=0$. Let
$\gamma:[a,b]\subset\Real\to M$ be a piecewise smooth curve on
$M$. The connection $\nabla$ defines the parallel transport
$\tau_\gamma:T_{\gamma(a)}M\to T_{\gamma(b)}M$, which is an
isomorphism of the pseudo-Euclidean spaces
$(T_{\gamma(a)}M,g_{\gamma(a)})$ and
$(T_{\gamma(b)}M,g_{\gamma(b)})$.

{\it The holonomy group} $G_x$ of $(M,g)$ at a point $x\in M$ is
the  Lie group that consists of the pseudo-orthogonal
transformations given by the parallel transports along all
piecewise smooth loops at the point $x$. It can be identified with
a Lie subgroup of the pseudo-orthogonal Lie group ${\rm
O}(r,s)={\rm O}(T_xM,g_x)$.

Recall that a tensor field $T$ on $(M,g)$ is {\it parallel} if $\nabla T=0$,
or equivalently $T$ is preserved by parallel transports: for any
piecewise smooth curve starting at a point $y\in M$ with an end-point
$z\in M$ it holds $\tau_\gamma T_y=T_z$, where $\tau_\gamma$ is the
extension of the parallel transport along $\gamma$ to tensors.

 The {\it fundamental principle for the holonomy groups}
 \cite{Besse}
states that there exists a one-to-one correspondents between
parallel tensor fields $T$ on $M$ and tensors $T_0$ of the same
type at $x$ preserved by the tensor extension of the
representation of the holonomy group.

\section{The de~Rham decomposition for Riemannian manifolds}\label{secdeR}

In this section we will give an algorithm allowing to decompose
(at least locally) any Riemannian manifold with not irreducible
holonomy group  in the product of a flat Riemannian manifold  and
of Riemannian manifolds with irreducible holonomy groups.  This
decomposition allows to restrict attention to Riemannian manifolds
with irreducible holonomy groups.

We consider a Riemannian manifold $(M,g)$ of dimension $n$ with
the holonomy group $G\subset{\rm O}(n)$ at a point $x\in M$.

If $(N,h)$ is another Riemannian manifold of dimension $m$ with
the holonomy group $H\subset{\rm O}(m)$ at a point $y\in N$, then
the product $(M\times N,g+h)$ is a Riemannian manifold with the
holonomy group $G\times H\subset {\rm O}(n+m)$ at the point
$(x,y)\in M\times N$; this holonomy group preserves the subspaces
$T_xM,T_yN\subset T_{(x,y)}M\times N$. This statement can be
inverted in the following way.

The {\it de~Rham decomposition Theorem} \cite{Besse} states that
if $(M,g)$ is simply connected and complete, and its holonomy
group is not irreducible, then $(M,g)$ can be decomposed into the
product of a flat Riemannian manifold $(M_0,g_0)$ and of
Riemannian manifolds $(M_1,g_1)$,...$(M_r,g_r)$ with irreducible
holonomy groups. For general $(M,g)$ with a not irreducible
connected holonomy group such decomposition exists only locally.
Thus the irreducibility of the connected holonomy group of a
Riemannian manifold $(M,g)$ is equivalent to the local
indecomposability of $(M,g)$.

Let us explain where the de~Rham decomposition comes from. Let
$x\in M$. Since the holonomy group $G\subset{\rm O}(n)$ is totally
reducible,  the tangent space $T_xM$ can be decomposed into an
orthogonal direct sum
\begin{equation}\label{decDeR}T_xM=E_{0x}\oplus E_{1x}\oplus\cdots\oplus
E_{rx},\end{equation} where $E_{0x}$ is the subspace consisting of
$G$-invariant vectors,  each subspace $E_{\alpha x}\subset T_xM$, $1\leq \alpha\leq r$,
is $G-$invariant  and the induced representation is irreducible.

The subspaces $E_{0 x},...,E_{r x}\subset T_xM$ define by means of
the parallel transport parallel distributions $E_0,...,E_r$ on
$M$. These distributions are involutive, and the manifolds
$M_0,...,M_r$ from the above decomposition are maximal integral
manifolds of these distributions passing through the point $x$.
The metrics $g_\alpha$ are the restrictions of $g$ to these
distributions.

The holonomy group of $(M,g)$ is the product
$$G=G_1\times\cdots\times G_r,$$ where $G_\alpha$ is the
restriction of $G$ to $E_{\alpha x}$. If the decomposition of
$(M,g)$ is global, then $G_\alpha$ is the holonomy group of the
manifold $(M_\alpha,g_\alpha)$. However, if the decomposition is
not global, then the holonomy group of $(M_\alpha,g_\alpha)$ may
be a proper subgroup of $G_\alpha$; in that case $G_\alpha$ is the
holonomy group of the induced connection on the distribution
$E_\alpha$ considered as a vector bundle over $M$.

The task is to find the distributions $E_\alpha$. We may find
$E_0$ as the distribution consisting of all parallel vector
fields. Then we may work with $E_0^\perp$ and $g$ restricted to
it. This allows us to assume that $E_0=0$.

 Note that if $(M,g)$ is indecomposable,  then the
dimension of parallel symmetric bilinear forms on $(M,g)$ equals
to one. This follows from the Fundamental principle for holonomy
groups and from the fact that any element in the second symmetric
power $\odot^2\Real^n$ of $T_xM\simeq \Real^n$ preserved by the
irreducible subgroup $G\subset{\rm O}(n)$ is proportional to the
metric $g_x$ at the point $x$.
 In general, the dimension of parallel symmetric bilinear forms on $(M,g)$
equals to $r$ (we assume that $E_0=0$) and this real vector space
is generated by $g_1,...,g_r$; here we assume that
$g_\alpha|_{E_\alpha\times E_\alpha}=g|_{E_\alpha\times E_\alpha}$
and $g_\alpha|_{E_\beta\times E}=0$ if $\alpha\neq \beta$. In
other words, \be\label{galfpr}
g_\alpha(\cdot,\cdot)=g(\pr_{E_\alpha}\cdot,\pr_{E_\alpha}\cdot),\ee
where the projection is taken with respect to decomposition
\eqref{decDeR}.

Finding all parallel symmetric bilinear forms on $(M,g)$ (e.g.
with Maple), we get an answer in the form $$c_1\tilde
g_1+\cdots+c_r\tilde g_r,$$ where $c_1,...,c_r\in\Real$ are
arbitrary and $\tilde g_1,...,\tilde g_r$ is a basis of the space
of all parallel symmetric bilinear forms on $(M,g)$. Since $g$ is
parallel, we may assume that $\tilde g_1=g$ (indeed, we may find a
linear independent subsystem in $\{g,\tilde g_1,...,\tilde g_r\}$
that contains $g$).

We may write
\begin{equation}\label{eqA}\tilde g_\alpha=\sum_{\beta=1}^rA_{\beta\alpha}
g_\beta,\qquad A_{\beta\alpha}\in\Real.\end{equation} Since we
need to find the numbers $A_{\beta\alpha}\in\Real$, we may work
with a fixed point $x\in M$. The proof of the following
proposition will allow to find the matrix $(A_{\beta\alpha})$
using some linear algebra.

\begin{prop} Let $V$ be a vector space with the Euclidean
metric $\eta$. Suppose that an orthogonal decomposition
\begin{equation}\label{decDeR1}V=V_{1}\oplus\cdots\oplus
V_{r}\end{equation} is fixed. Let
$$\eta_\alpha(\cdot,\cdot)=\eta(\pr_{V_\alpha}\cdot,\pr_{V_\alpha}\cdot),\quad
\alpha=1,...,r.$$  If an arbitrary basis
$\tilde\eta_1,...,\tilde\eta_r$ of the vector space
$\spa\{\eta_1,...,\eta_r\}$ is given, then the forms $\eta_\alpha$
and the decomposition \eqref{decDeR1} can be reconstructed up to a
permutation.
\end{prop}

{\bf Proof.} We obtain the relation
$$\tilde\eta_\alpha=\sum_{\beta=1}^rA_{\beta\alpha}
\eta_\beta,\qquad A_{\beta\alpha}\in\Real,\quad \alpha=1,...,r.$$
We may assume that $\tilde\eta_1=\eta$.
 Consider $\tilde \eta_{2}$, then
$$V=F\oplus F^\bot,$$ where $$F=\ker\tilde\eta_2=\{X\in V|\tilde \eta_{2}(X,Y)=0\text{ for all } Y\in V\}$$
is the kernel of $\tilde \eta_{2}$ and $F^\bot $ is its orthogonal
complement with respect to $\eta$. Both $F$ and $F^\bot $ consist
of some of $V_{\alpha}$, i.e. this decomposition is orthogonal
with respect to all tensors $\eta_{\alpha}$. Consider the
decomposition $V=F\oplus F^\bot$, take the restrictions of $\tilde
\eta_{3}$ to each of these spaces and decompose $F$ and $F^\bot $
in the same manner. Continue this process for all $\tilde
\eta_{\alpha}$, then we get a decomposition
\begin{equation}\label{decF}V=F_1\oplus\cdots\oplus F_s\end{equation} such that the restriction of
each $\tilde \eta_{\alpha}$ to any of $F_k$ is either zero or
non-degenerate. Now we continue to subdivide this decomposition.
Let $\alpha$ run from $2$ to $r$.  Consider the restrictions
$\tilde \eta_{\alpha}|_{F_k\times F_k}$ and $\eta|_{F_k\times
F_k}$ of $\tilde \eta_{\alpha}$ and $\eta$ to each $F_k$. If the
restriction $\tilde \eta_{\alpha}|_{F_k\times F_k}$  is non-zero
(i.e. $\tilde \eta_{\alpha}|_{F_k\times F_k}$ is non-degenerate),
and $\tilde \eta_{\alpha}|_{F_k\times F_k}$ is not proportional to
$\eta|_{F_k\times F_k}$, then we consider a new
$\tilde\eta_\alpha$: its restriction to $F_l\times F_l$ for $l\neq
k$ remains the same, and we change its restriction $\tilde
\eta_{\alpha}|_{F_k\times F_k}$ to $ \tilde
\eta_{\alpha}|_{F_k\times F_k}- b \eta|_{F_k\times F_k}$, where
$b\in \Real$ is a number such that the restriction  $(\tilde
\eta_{\alpha }-b \eta)|_{F_k\times F_k}$ to $F_k$ is degenerate.
To find such $b$, take any vector $X\in F_k$ such that
$\eta(X,X),\tilde \eta_{\alpha}(X,X)\neq 0$ (in the case of
positive definite $\eta$, any non-zero $X$ satisfies this
condition) and set $b=\frac{\tilde
\eta_{\alpha}(X,X)}{\eta(X,X)}$.
 Using the new tensor $\tilde \eta_{\alpha}$, we
may subdivide $F_k$ (since now $\tilde \eta_{\alpha}|_{F_k\times
F_k}$ is degenerate and non-zero), i.e. we subdivide  the
decomposition \eqref{decF}. Continue this process. At the end we
will get that  the restriction of any $\tilde \eta_{\alpha}$ to
any  $F_k$ is either zero, or  it is proportional to
$\eta|_{F_k\times F_k}$. This means that the number $s$ in
decomposition \eqref{decF} equals $r$, that is, decomposition
\eqref{decF} is the decomposition \eqref{decDeR1},
$F_\alpha=V_{\alpha}$ (up to a renumbering). Then we may find the
forms $\eta_\alpha$, and using the initial forms
$\tilde\eta_\alpha$, we find the matrix $(A_{\beta\alpha})$.
Finally, $V_\alpha=(\ker \eta_{\alpha})^{\bot_\eta}$. This proves
the proposition. $\Box$

Thus using \eqref{eqA} considered at the point $x$, we find the
matrix $(A_{\beta \alpha})$. Then using \eqref{eqA} and the
inverse matrix, we find the metrics $g_\alpha$. Now for any $y\in
M$ find $E_{\alpha y}$ using one of the formulas:
\be\label{Ealfy1}E_{\alpha y}=(\ker g_{\alpha y})^{\bot_g},\quad
\text{or}\ee \be\label{Ealfy}E_{\alpha
y}=\cap_{\beta\neq\alpha}\ker g_{\beta y}=\{X\in T_yM|g_{\beta
y}(X,\cdot)=0 \, \text{ for all } \beta\neq \alpha\}.\ee Thus we
know the distributions $E_\alpha$.

\subsection{Example} In order to find a decomposable metric, we simply take the local
metric on the product of two spheres:
$$g=(dy_1)^2+\sin^2y_1(dy_2)^2+(dy_3)^2+\sin^2y_3(dy_4)^2,\quad 0\leq y_1,y_3\leq \pi.$$
Consider the new coordinates $$x_1=y_1-y_3,\quad x_2=y_2-y_4,\quad
x_3=y_3,\quad x_4=y_4.$$ The metric takes the form
$$g=\left(\begin{array}{cccc}
1 & 0 & 1 & 0\\
0 & \sin^2(x_1+x_3) & 0 & \sin^2(x_1+x_3)\\
1 & 0 & 2 & 0\\
0 & \sin^2(x_1+x_3) & 0 & \sin^2 x_3 +\sin^2(x_1+x_3)
\end{array}\right).$$
Looking at this metric, it is not obvious that it is decomposable.
Now we will apply the above algorithm in order to decompose the
obtained metric. Using Maple, we get that there are no non-zero
parallel vector fields and parallel symmetric bilinear forms are
the following:
$$\left(\begin{array}{cccc}
C_2 & 0 & C_2 & 0\\
0 & C_2\sin^2(x_1+x_3) & 0 & C_2\sin^2(x_1+x_3)\\
C_2 & 0 & C_1 & 0\\
0 & C_2\sin^2(x_1+x_3) & 0 & C_1\sin^2 x_3
+C_2(\sin^2(x_1+x_3)-\sin^2x_3)
\end{array}\right).$$ Let $\tilde g_1=g,\quad \tilde g_2=
 \left(\begin{array}{cccc}
0 & 0 & 0 & 0\\
0 & 0 & 0 & 0\\
0 & 0 & 1 & 0\\
0 & 0 & 0 & \sin^2 x_3
\end{array}\right).$
Consider the point $x=\left(0,0,\frac{\pi}{2},0\right)$. We obtain
the decomposition $$T_xM=F\oplus F^\perp=\ker \tilde g_{2x}\oplus
(\ker \tilde g_{2x})^{\perp},$$ $$ \ker \tilde
g_{2x}=\spa\{(\partial_{x_1})_x,(\partial_{x_2})_x\},\quad (\ker
\tilde
g_{2x})^{\perp}=\spa\{(\partial_{x_1}-\partial_{x_3})_x,(\partial_{x_2}-\partial_{x_4})_x\}.$$
Since the number of the summands in this decomposition equals to
the dimension of parallel symmetric bilinear forms, we get
$$E_{1x}=\ker \tilde g_{2x},\quad E_{2x}=(\ker \tilde
g_{2x})^{\perp}.$$ Using \eqref{galfpr} evaluated at $x$, we
obtain
$$g_{1x}= \left(\begin{array}{cccc}
1 & 0 & 1 & 0\\
0 & 1 & 0 & 1\\
1 & 0 & 1 & 0\\
0 & 1 & 0 & 1
\end{array}\right),\quad g_{2x}=\left(\begin{array}{cccc}
0 & 0 & 0 & 0\\
0 & 0 & 0 & 0\\
0 & 0 & 1 & 0\\
0 & 0 & 0 & 1
\end{array}\right).$$
We conclude that $g_x=g_{1x}+g_{2x}$, $\tilde g_{2x}=g_{2x}$, and
$$g_1=g-\tilde g_2,\quad g_2=\tilde g_2.$$ Using \eqref{Ealfy}, we get that the
distribution $E_1$ is spanned by the vector fields
$\partial_{x_1}$, $\partial_{x_2}$, and $E_2$ is spanned by the
vector fields $\partial_{x_1}-\partial_{x_3}$,
$\partial_{x_2}-\partial_{x_4}$. We find the coordinates
$z_1,...,z_4$ adopted to the decomposition $TM=E_1\oplus E_2$
requiring that
$$\partial_{z_1}=\partial_{x_1},\quad
\partial_{z_2}=\partial_{x_2},\quad
\partial_{z_3}=\partial_{x_1}-\partial_{x_3},\quad
\partial_{z_4}=\partial_{x_2}-\partial_{x_4},$$ and we get
$$z_1=x_1+x_3,\quad z_2=x_2+x_4,\quad z_3=-x_3,\quad z_4=-x_4.$$
With respect to these coordinates it holds
$$g=(dz_1)^2+\sin^2z_1(dz_2)^2+(dz_3)^2+\sin^2 z_3(dz_4)^2,$$
i.e. we obtain the initial metric.

\section{The Wu decomposition for Lorentzian manifolds}\label{secWu}

{\it The Wu decomposition Theorem} \cite{Wu} generalizes the
de~Rham Theorem for the case of pseudo-Riemannian manifolds. It
states that a pseudo-Riemannian manifold $(M,g)$ with not weakly
irreducible connected holonomy group can be decomposed at list
locally into the product of pseudo-Riemannian manifolds
$(M_0,g_0),$...,$(M_r,g_r)$ such that $(M_0,g_0)$ is flat and the
holonomy groups of $(M_1,g_1),$...,$(M_r,g_r)$ are weakly
irreducible. Recall that a subgroup of a pseudo-orthogonal group
is weakly irreducible, if it does not preserve any proper
non-degenerate subspace of the pseudo-Euclidean space.
Consequently, a pseudo-Riemannian manifold $(M,g)$ is locally
indecomposable if and only if its connected holonomy group is
weakly irreducible. In that case $E_0$ is not the distribution
that consists of all parallel vector fields on $(M,g)$, but it is
a non-degenerate subdistribution of the last one and in general it
is not defined uniquely. Note that if the holonomy group is weakly
irreducible and not irreducible, then it preserves an isotropic
subspace of the tangent space, but this does not imply the local
decomposability of the manifold.

The algorithm of Section \ref{secdeR} works also for
pseudo-Riemannian manifolds if we know that the holonomy group of
each factor in the decomposition is irreducible. In that case the
dimension of parallel symmetric bilinear forms equals to the
number of the manifolds in the decomposition  (without loss of
generality we assume that $E_0=0$), then \eqref{eqA} holds. The
problem is that a locally indecomposable pseudo-Riemannian
manifold may admit  parallel symmetric bilinear forms not
proportional to the metric. The structure of such forms is known
only in the case of Lorentzian manifolds.

Now we consider a Lorentzian manifold $(M,g)$. And we will obtain
the Wu decomposition of that manifold.

A.V.~Aminova \cite{Aminova} proved that if a locally
indecomposable Lorentzian manifold $(M,g)$ admits a parallel
symmetric bilinear form not proportional to the metric, then
$(M,g)$ admits a parallel light-like vector field $p$,  the space
of parallel bilinear forms is 2-dimensional, and it is spanned by
$g$ and $\theta\otimes\theta$, where $\theta=g(p,\cdot)$ is the
1-form corresponding to $p$ (see also the paper \cite{H88} by
G.S.~Hall for the case of dimension 4).

It is clear that in the Wu decomposition of a Lorentzian manifold
only one submanifold is Lorentzian, and all the other are
Riemannian. The Lorentzian part is locally indecomposable and
admits a parallel light-like vector field if and only if the
restriction of $g$ to the space of parallel vector fields is
degenerate (this property may be checked at a single point if we
restrict the parallel vector fields to that point). If the
Lorentzian part is contained in the distribution $E_0$, or it does
not admit a parallel light-like vector field, then the algorithm
of Section \ref{secdeR} works.

Suppose that the Lorentzian part (that we assume to be
$(M_r,g_r)$) in the Wu decomposition is indecomposable and admits
a parallel vector field $p$. Let $\tilde E_0\subset TM$ be the
subbundle spanned by all parallel vector fields on $(M,g)$. Then,
$p\in\Gamma(\tilde E_0)$. Let $E_{0x}\subset\tilde E_{0x}$ be any
subspace complementary to $\Real p_x$. Let $E_0\subset\tilde E_0$
be the subbundle spanned by parallel vector fields with values in
$E_{0x}$ at the point $x$. Then $E_0\subset TM$ is a parallel
subbundle and the restriction of $g$ to $E_0$ is non-degenerate.
We consider $E_0^\perp$ and the restriction of $g$ to it. Hence we
again may assume that $E_0=0$. Then the space of parallel bilinear
forms on $(M,g)$ is spanned by $g_1,...,g_r,\theta\otimes\theta$
and it is of dimension $r+1$.

Finding all parallel symmetric bilinear forms on $(M,g)$, we get
an answer in the form $$c_1\bar g_1+\cdots+c_{r+1}\bar g_{r+1},$$
where $c_1,...,c_{r+1}\in\Real$  are arbitrary and $\bar
g_1,...,\bar g_{r+1}$ is a basis of the space of all parallel
symmetric bilinear forms on $(M,g)$. Since $g$ is parallel, we may
assume that $\bar g_1=g$.

 There exist real numbers
$(C_{\beta\alpha})_{\beta,\alpha=1}^{r+1}$ such that
\be\label{eqbarg}\bar
g_{\alpha}=\sum_{\beta=1}^{r}C_{\beta\alpha}g_\beta+C_{r+1\,\alpha}\theta\otimes\theta.\ee

Again, we may find the matrix $C_{\beta\alpha}$ sitting at a
single point and doing linear algebra.

\begin{prop} Let $V$ be a vector space with the Minkowski
metric $\eta$. Suppose that an orthogonal decomposition
\begin{equation}\label{decWu1}V=V_{1}\oplus\cdots\oplus
V_{r}\end{equation} is fixed such that the restriction of $\eta$
to $V_r$ is of Minkowski signature. Let $p\in V_r$ be a fixed
non-zero isotropic vector, $\theta=g(p,\cdot)$, and let
$$\eta_\alpha(\cdot,\cdot)=\eta(\pr_{V_\alpha}\cdot,\pr_{V_\alpha}\cdot),\quad
\alpha=1,...,r.$$ If an arbitrary basis
$\bar\eta_1,...,\bar\eta_{r+1}$ of the vector space
$\spa\{\eta_1,...,\eta_r,\theta\otimes\theta\}$ is given, then the
forms $\eta_\alpha$ and the decomposition \eqref{decWu1} can be
reconstructed up to a permutation of the subspaces
$V_1,...,V_{r-1}$.
\end{prop}

{\bf Proof.} Again we may assume that $\bar\eta_1=\eta$. Consider
the relation \be\label{eqbarg1}\bar
\eta_{\alpha}=\sum_{\beta=1}^{r}C_{\beta\alpha}\eta_\beta+C_{r+1\,\alpha}\theta\otimes\theta.\ee
In particular, $(C_{\beta 1})_{\beta=1}^{r+1}=(1,0,...,0).$

Let $q\in V$ be a light-like vector not proportional to $p$, i.e.
$\theta(q)=\eta(p,q)\neq 0$. To find such vector it is enough to
take any vector $X\in V$ such that $\eta(p,X)\neq 0$, and if
$\eta(X,X)\neq 0$, then take $q=p-\frac{2\eta(p,X)}{\eta(X,X)}X$.
It is clear that
$$\bar \eta_{\alpha
}(p,q)=C_{r\alpha}\eta_{r}(p,q)=C_{r\alpha}\eta(p,q),\quad 2\leq
\alpha\leq r+1.$$ This allows to find the coefficients
$C_{r\alpha}$. Changing each $\bar\eta_\alpha$ to
$\bar\eta_\alpha-C_{r\alpha}\eta$, we obtain that $C_{r\alpha}=0$
for $2\leq \alpha\leq r+1$.

The following three  lemmas will allow to get the algorithm.

\begin{lem} Let $2\leq \alpha\leq r+1$. If $C_{r+1\,\alpha} =0$, then
$$\ker \bar \eta_{\alpha }=\oplus_{1\leq \beta\leq r,\, C_{\beta\alpha}=0} V_{\beta },$$
and $$\ker \eta|_{\ker \bar \eta_{\alpha }\times\ker \bar
\eta_{\alpha }}=0.$$ If  $C_{r+1\,\alpha}\neq 0$, then
$$\ker \bar \eta_{\alpha }=\oplus_{1\leq \beta\leq r-1,\, C_{\beta\alpha}=0} V_{\beta }\oplus\{X\in V_{r}|\theta(X)=0\},$$
and $$\ker \eta|_{\ker \bar \eta_{\alpha }\times\ker \bar
\eta_{\alpha }}=\Real p.$$\end{lem} {\it Proof.} Suppose that
$C_{r+1\,\alpha} =0$. If $C_{\beta\alpha}= 0$ then it is clear
that $V_{\beta }\subset\ker\bar \eta_{\alpha }$. Let $X\in\ker\bar
\eta_{\alpha }$. We may write $X=X_1+\cdots+X_r$, where
$X_\gamma\in V_{\gamma }$. Suppose that $C_{\beta\alpha}\neq 0$.
Let $Y\in V_{\beta }$ be any vector. Then, $0=\bar \eta_{\alpha
}(X,Y)=C_{\beta \alpha}\eta_\beta(X_\beta,Y)$ for any $Y\in
V_{\beta }$, i.e. $X_\beta=0$. This implies the first equality of
the first statement. The second equality is obvious, since $\eta$
is non-degenerate on each $V_{\beta }$.

Suppose that $C_{r+1\,\alpha}\neq 0$. The inclusion $\supset$ in
the first equality  is obvious. Let $X\in\ker\bar \eta_{\alpha }$.
We write $X=X_1+\cdots+X_r$, where $X_\gamma\in V_{\gamma }$.
Suppose that $C_{\beta\alpha}\neq 0$, $1\leq\beta\leq r-1$. As
above, this implies $X_\beta=0$. Let $Y\in V_{r}$ be a vector such
that $\theta(Y)\neq 0$. Then $0=\bar
\eta_\alpha(X,Y)=C_{r+1\,\alpha}\theta(X_r)\theta(Y)$, i.e.
$\theta(X_r)=0$. This proves the second statement. \hfill $\Box$

This lemma allows us easily indicate whether  $C_{r+1\,\alpha} =0$
or not: we should compute $\ker \bar \eta_{\alpha }$ and $\ker
\eta|_{\ker \bar \eta_{\alpha }\times\ker \bar \eta_{\alpha }}$.
If $C_{r+1\,\alpha} =0$ for some $\alpha\geq 2$, then we add
$\theta\otimes \theta$ to $\bar \eta_\alpha$. Then by the lemma we
get
$$\ker \bar \eta_{\alpha }=\oplus_{1\leq \beta\leq r-1,\, C_{\beta\alpha}=0} V_{\beta }\oplus\{X\in V_{r}|\theta(X)=0\},\quad 2\leq \alpha\leq r+1.$$

Let us consider the following vector space:
$$W=\cap_{\alpha=2}^{r+1} \ker\bar \eta_{\alpha }.$$

\begin{lem}
It   holds $$W=\{X\in V_{r}|\theta(X)=0\},$$
$$W^{\perp_\eta}=V_{1}\oplus\cdots\oplus V_{r-1\, }\oplus\Real p.$$
\end{lem}

{\it Proof.} We claim that for any $\beta$, $1\leq \beta\leq r-1$
there exists an $\alpha$, $2\leq \alpha\leq r+1$, such that
$C_{\beta\,\alpha}\neq 0$. Indeed, if the claim is wrong then
there exists a $\beta$, $1\leq \beta\leq r-1$ such that it holds
$$(C_{\beta \alpha})_{\alpha=1}^{r+1}=(C_{r
\alpha})_{\alpha=1}^{r+1}=(1,0,...,0),$$ i.e. the matrix
$(C_{\beta \alpha})_{\beta,\alpha=1}^{r+1}$ is degenerate, that
gives a contradiction. The first equality follows from the claim.
The second equality is obvious. \hfill $\Box$

\begin{lem}\label{lemalg3} The intersection $\cap_{\alpha=2}^{r+1}(W^{\perp_\eta})^{\perp_{\bar \eta_\alpha}}$ that can be written as
$$\{X\in V|
\bar \eta_{\alpha }(X,Y)=0\text{ for all } 2\leq \alpha\leq
r+1,\,Y\in W^{\perp_\eta}\}$$ coincides with $V_{r}$. \end{lem}
{\it Proof.}  Suppose that $X\in
\cap_{\alpha=2}^{r+1}(W^{\perp_\eta})^{\perp_{\bar \eta_\alpha}}$.
Then for any $\alpha$,  $2\leq\alpha\leq r+1$ and all $Y\in
W^{\perp_\eta}$ it holds $\bar \eta_{\alpha }(X,Y)=0$. Consider
the decomposition $X=X_1+\cdots+X_r$. Let $1\leq\beta\leq r-1$.
Above we have seen that there exist $\alpha$, $2\leq\alpha\leq
r+1$ such that $C_{\beta\alpha}\neq 0$. Let $Y\in V_{\beta
}\subset W^{\perp_\eta}$. Then,
$$0=\eta_{\alpha }(X,Y)=C_{\beta\alpha}\eta_{\beta }(X,Y)=C_{\beta\alpha}\eta_{\beta }(X_\beta,Y).$$
Consequently, $X_\beta=0$, and $X=X_r\in V_{r}$.

Conversely, if $X\in V_{r}$, $2\leq \alpha\leq r+1$, $Y\in
W^{\perp_\eta}$, then $\eta_{\alpha
}(X,Y)=C_{r+1\,\alpha}\theta(X)\theta(Y)=0$, since $\theta(Y)=0$.
This proves the lemma.  \hfill $\Box$

Thus in order to find the space  $V_{r}$, it is enough to compute
the spaces
 $$W=\cap_{\alpha=2}^{r+1} \ker\bar \eta_{\alpha },\quad W^{\perp_\eta},\quad \cap_{\alpha=2}^{r+1}(W^{\perp_\eta})^{\perp_{\bar \eta_\alpha}}.$$

 Now we consider $(V_{r})^{\perp_\eta}$ and the restrictions of the forms $\bar \eta_{1}=\eta, \bar \eta_{2},...,\bar \eta_{r+1 \, }$ to  $(V_{r})^{\perp_\eta}$. Clearly, the rank of this system equals to $r-1$.
  Let $\tilde \eta_{1},..., \tilde \eta_{r-1\,}$ be a linearly independent subsystem such that $\tilde \eta_{1}=\eta$.
 There exist real numbers  $(A_{\beta\alpha})_{\beta,\alpha=1}^{r-1}$ such that
\begin{equation}\label{eqA1}\tilde \eta_{\alpha }=\sum_{\beta=1}^{r-1}A_{\beta\alpha}
\eta_{\beta }.\end{equation} These numbers can be found in the
same way as in Section \ref{secdeR}. Using \eqref{eqA1} and the
inverse matrix to $A_{\beta\alpha}$,  we find the forms $\eta_{1
},...,\eta_{r-1\, }$. Next,
$\eta_{r}=\eta-\eta_{1}-\cdots-\eta_{r-1\,}$ and
$\theta=\eta(p,\cdot)=\eta_{r}(p,\cdot)$. Using the initial
$\bar\eta_1,...,\bar\eta_{r+1}$ and \eqref{eqbarg1}, we find the
matrix $(C_{\beta\alpha})$. Now, $V_\alpha=(\ker
\eta_{\alpha})^{\bot_\eta}$. The proposition is proved. $\Box$

Thus we have found the matrix $(C_{\beta\alpha})$ from the
equation \eqref{eqbarg}.
 Then using \eqref{eqbarg} and the
inverse matrix to $(C_{\beta\alpha})$, we find the metrics
$g_\alpha$. Now the distributions $V_\alpha$ can be found using
\eqref{Ealfy1} or \eqref{Ealfy}.

\subsection{Example} Consider the following product metric
$$g=(dy_1)^2+\sin^2 y_1(dy_2)^2+2dy_3dy_5+(dy_4)^2+y_4^2(dy_5)^2.$$
After the coordinate transformation $$x_1=y_1-y_4,\quad
x_2=y_2,\quad ...,\quad x_5=y_5$$ the metric takes the form
$$
g=\left(\begin{array}{ccccc}
1 & 0 & 0 & 1 & 0 \\
0 & \sin^2(x_1+x_4) & 0 & 0 & 0 \\
0 & 0 & 0 & 0 & 1 \\
1 & 0 & 0 & 2 & 0 \\
0 & 0 & 1 & 0 & x_4^2
\end{array}\right).$$
Now we apply the algorithm of this section to the obtained metric.
The space of parallel vector fields is one-dimensional and it is
spanned by the light-like vector field  $p=\partial_{x_3}$. Hence,
$\theta=dx_5$. The space of parallel bilinear symmetric forms is
three-dimensional and consists of the following matrices:
$$
\left(\begin{array}{ccccc}
C_1 & 0 & 0 & C_1 & 0 \\
0 & C_1\sin^2(x_1+x_4) & 0 & 0 & 0 \\
0 & 0 & 0 & 0 & C_3 \\
C_1 & 0 & 0 & C_1+C_3 & 0 \\
0 & 0 & C_3 & 0 & C_2+C_3x_4^2
\end{array}\right).$$
It is clear that in the above notation $r=2$. Let $ \bar
g_1=g,\quad \bar g_2=\left(\begin{array}{ccccc}
0 & 0 & 0 & 0 & 0 \\
0 & 0 & 0 & 0 & 0 \\
0 & 0 & 0 & 0 & 0 \\
0 & 0 & 0 & 0 & 0 \\
0 & 0 & 0 & 0 & 1 \end{array}\right),\quad
\bar g_3=
 \left(\begin{array}{ccccc}
0 & 0 & 0 & 0 & 0 \\
0 & 0 & 0 & 0 & 0 \\
0 & 0 & 0 & 0 & 1 \\
0 & 0 & 0 & 1 & 0 \\
0 & 0 & 1 & 0 & x^2_4
\end{array}\right).$
We obtain that $\bar g_2=\theta\otimes\theta$. Consider the
decomposition $$\bar
g_3=C_{13}g_1+C_{23}g_2+C_{33}\theta\otimes\theta.$$ Let
$x=\left(\frac{\pi}{2},0,0,0,0\right)$. Note that the vector
$q_x=(\partial_{x_5})_x$ is light-like and satisfies
$g_x(p_x,q_x)=1$. Since $\bar g_{3x}(p_x,q_x)=1$, we get
$C_{23}=1$. We change $\bar g_3$ to $\bar g_3-g$. Then $C_{23}=0$.
Next, $$\ker \bar
g_{3x}=\spa\{(\partial_{x_1}-\partial_{x_4})_x,(\partial_{x_3})_x,(\partial_{x_5})_x\},\quad
\ker (g_x|_{\ker \bar g_{3x}\times\ker \bar g_{3x}})=0,$$ i.e.
$C_{33}=0$. We change  $\bar g_3$ to $\bar
g_3+\theta\otimes\theta$, then $C_{33}=1$ and
$$ \bar g_3=\left(\begin{array}{ccccc}
-1 & 0 & 0 & -1 & 0 \\
0 & -\sin^2(x_1+x_3) & 0 & 0 & 0 \\
0 & 0 & 0 & 0 & 0 \\
-1 & 0 & 0 & -1 & 0 \\
0 & 0 & 0 & 0 & 1
\end{array}\right).$$
Now, $$W=\ker \bar g_{2x}\cap \ker \bar
g_{3x}=\spa\{(\partial_{x_1}-\partial_{x_4})_x,(\partial_{x_3})_x\},\qquad
W^{\perp_g}=\spa\{(\partial_{x_1})_x,(\partial_{x_2})_x,(\partial_{x_3})_x\}.$$
Using Lemma \ref{lemalg3}, we find
$$E_{2x}=\spa\{(\partial_{x_1}-\partial_{x_4})_x,(\partial_{x_3})_x,(\partial_{x_5})_x\}.$$
Consequently,
$E_{1x}=E_{2x}^\perp=\spa\{(\partial_{x_1})_x,(\partial_{x_2})_x\}.$
Using \eqref{galfpr} evaluated at the point $x$, we get $$
g_{1x}=\left(\begin{array}{ccccc}
1 & 0 & 0 & 1 & 0 \\
0 & 1 & 0 & 0 & 0 \\
0 & 0 & 0 & 0 & 0 \\
1 & 0 & 0 & 1 & 0 \\
0 & 0 & 0 & 0 & 0
\end{array}\right),\quad g_{2x}=
\left(\begin{array}{ccccc}
0 & 0 & 0 & 0 & 0 \\
0 & 0 & 0 & 0 & 0 \\
0 & 0 & 0 & 0 & 1 \\
0 & 0 & 0 & 1 & 0 \\
0 & 0 & 1 & 0 & 0
\end{array}\right).$$
Hence, $\bar g_{3x}=-g_{1x}+\theta_x\otimes\theta_x$. That implies
$$
g_{1}=-\bar g_3+\theta\otimes\theta=\left(\begin{array}{ccccc}
1 & 0 & 0 & 1 & 0 \\
0 & \sin^2(x_1+x_4) & 0 & 0 & 0 \\
0 & 0 & 0 & 0 & 0 \\
1 & 0 & 0 & 1 & 0 \\
0 & 0 & 0 & 0 & 0
\end{array}\right),\quad g_{2}=g-g_1=
\left(\begin{array}{ccccc}
0 & 0 & 0 & 0 & 0 \\
0 & 0 & 0 & 0 & 0 \\
0 & 0 & 0 & 0 & 1 \\
0 & 0 & 0 & 1 & 0 \\
0 & 0 & 1 & 0 & x^2_4
\end{array}\right).$$
Using \eqref{Ealfy}, we see that the distribution $E_1$ is spanned
by the vector fields $\partial_{x_1}$,  $\partial_{x_2}$, and the
distribution $E_2$ is spanned by the vector fields
$\partial_{x_1}-\partial_{x_4}$,  $\partial_{x_3}$,
$\partial_{x_5}$. We find the coordinates $z_1,...,z_5$ adopted to
the decomposition $TM=E_1\oplus E_2$ requiring that
$$\partial_{z_1}=\partial_{x_1},\quad
\partial_{z_2}=\partial_{x_2},\quad
\partial_{z_3}=\partial_{x_3},\quad
\partial_{z_4}=\partial_{x_1}-\partial_{x_4},\quad \partial_{z_5}=\partial_{x_5},$$ and we get
$$z_1=x_1+x_4,\quad z_2=x_2,\quad ...,\quad  z_5=x_5.$$
With respect to these coordinates it holds
$$g=(dz_1)^2+\sin^2z_1(dz_2)^2+2dz_3dz_5+(dz_4)^2+z_4^2(dz_5)^2,$$
and we obtain the initial metric.

\end{document}